# THE HISTORY OF THE TOTAL CHROMATIC NUMBER CONJECTURE

Hossein Shahmohamad
School of Mathematical Sciences
Rochester Institute of Technology, Rochester, NY 14623
E-mail: **hxssma@rit.edu**

## Abstract

The total chromatic number conjecture which has appeared in a few hundred articles and in numerous books thus far is now one of the classic mathematical unsolved problems. It appears that many authors coincidentally have attributed it to Professor M. Behzad and/or to Professor V. G. Vizing. Eventually after four decades, Professor A. Soifer investigated the origin of this conjecture; published his findings in The Mathematical Coloring Book - 2009; and stated that, "In my opinion this unquestionably merits the joint credit to Vizing and Behzad." After checking all the arguments presented and the blames cited, I decided to investigate the controversy stated in this book on my own. My findings which are presented in this report specifically signify the following two points.

- M. Behzad is the sole author of the Total Chromatic Number Conjecture.

- The wrong referrals provided by numerous authors over the last forty four years, to indicate Vizing's authorship, must be brought to the attention of the authors and researchers, by appropriate means, as soon as possible.

## Introduction

The Total Chromatic Number Conjecture is stated as follows: Color all the vertexes and all the edges of a given graph G simultaneously in such a way that no two adjacent vertexes have the same color, no two adjacent edges have the same color, and that the color of each edge is different from the colors of its end vertexes. Then the minimum number of colors required to do so is at most the maximum degree of G plus two.

Aside from Abstract and Introduction, this report contains the following sections. The Inquiry: consists of the main part of the e-mail I sent to Professor Vizing on December 19, 2010. In this self explanatory section the subject is introduced. The necessary references presented at the end of this report were originally part of this Inquiry. The Reply: contains Vizing's reply, sent on December 31, 2010, is quite short and decisively convincing. The Follow-up: consists of a copy of the unanswered e-mail sent on January 11, 2011. Concluding Remarks: is the last section of the report.

# *The Inquiry*

From: Hossein Shahmohamad
Sent: Sunday, December 19, 2010 3:20 PM
To: vizing@paco.net
Subject: Total Chromatic Number Conjecture

Dear Professor Vizing:

 I wish you a very happy New Year 2011. Hopefully, for many years to come, you will continue producing great results in mathematics.

 After I read Section 16.2 of The Mathematical Coloring Book, written by Professor A. Soifer (A.S.) and noticed the accusations therein, I decided to investigate the origin of the Total Chromatic Number Conjecture, TCNC, on my own. As you well know, TCNC has appeared in several books and in a few hundred articles thus far. Authors have attributed this conjecture to Professor Mehdi Behzad (M.B.) and/or to you. The main references are the Ph.D. thesis of M.B. [1] and the 1964, 1965, or 1968 articles [5,6,7] written by you. Out of these three articles, the first two [5,6] clearly do not contain even a single word about TCNC! Although several distinguished graph theorists, including obviously yourself, are well aware of this fact, to my best knowledge, no one has made any attempt to rectify these commonly wrong referrals. Isn't it one of our main professional as well as academic responsibilities to prevent such mistakes?

Now, the situation about your 1968 paper [7] is mind boggling indeed. As mentioned in its introduction, not all the unsolved problems therein are posed by the author. In fact, one of the references for the unsolved problems introduced therein is the Proceedings of the International Symposium - Rome which was held in 1966, and it contains [2] referring to the original work of M.B. [1]. In contrast to several other problems presented in the paper [7] you did not specify an author for the TCNC mentioned in Section 5. Thus, unfortunately, one may think that TCNC has been copied from reference [2]. This ambiguity may be one of the main reasons inducing demand for clarification by the anonymous referee mentioned in Section 16.2 of A.S.'s book [4].

Another point is that A. A. Zykov (A.A.Z.) attributed TCNC, as an unsolved problem, to you [8] and provided the wrong 1965 reference [6]. This occurred at the Maneback Colloquium held in 1967. Please note that at the end of your 1964 article [5] A.A.Z. was thanked for his valuable comments. So A.A.Z. must have been aware of the details of your work. However, he committed such a mistake; i.e., referred to a paper that does not contain a single word on TCNC.

Now it is interesting to note that A.A.Z. also attended the International Symposium, Rome held in 1966 and presented two papers, one of which was: Some New Results of Soviet mathematicians. In this article [9] he referred to four of your work all different from TCNC. The real big question mark here is: what did happen between 1966 and 1967 that TCNC deserved to be presented in Maneback, but not in Rome. Moreover, don't you think that A.A.Z. should have given credit for TCNC to M.B. as well? Obviously, the same question applies to you, Professor Vizing, in regards to your 1968 article.

From your paper [7] I quote," To conclude this section we mention a form of the coloring of multigraphs that has not been investigated."  Therefore, total coloring was not investigated then, and as far as I have researched, I have not been able to find any result from you on total concepts published elsewhere. Then, clearly one wonders how the conjecture was posed.

Please note that according to your celebrated theorem [5] the edge chromatic number of a graph is at most one more than its maximum degree. It is not easy to believe that if one imposes two additional restricting conditions on the simultaneous coloring of the edges and the vertexes, then the relevant parameter will always be at most two more than its maximum degree. So, it would be hard to imagine that your theorem about edge coloring has inspired the introduction of the TCNC.

Dear Professor Vizing: What has been described are just a few out of many points related to the TCNC history. But, to summarize, in my opinion, the way TCNC has been presented by several colleagues, including A.A.Z, A.S. and yourself has raised serious questions and skepticism. Not only the wrong initial referral provided by A.A.Z. in 1967 has been continuously ignored till now, I surprisingly found a similar overlook in your own interview published in 2000 [3].

As you have noticed, I am not writing this e-mail to accuse anybody by any means. Instead, I am simply reporting my findings, indicating that TCNC has been mishandled by the authorities since its inception. In particular, I would like to add that, in spite of all aforementioned negligence, the title, "Total Insanity around Total Chromatic Number Conjecture" still seems to be an inappropriate and totally unacceptable choice for Section 16.2 of The Mathematical Coloring Book [4].

Your comments will be highly appreciated. If I can be of any help please let me know. I wish you all the best.

 Sincerely Yours,

Hossein Shahmohamad

## ***The Reply***

From: V. Vizing [mailto:vizing@paco.net]
Friday, December 31, 2010 4:35 AM
To: Hossein Shahmohamad
Subject: Re: Total Chromatic Number Conjecture

Dear Professor Shahmohamad,

Happy New Year to you! Concerning total chromatic number:

I published about it only in "Some Unsolved Problems..." (1968) and "Interview..." (2000). So other references are wrong. I didn't read the book of Soifer. All the best in 2011!

Sincerely Yours,

    V. Vizing

## *The Follow-up*

From: Hossein Shahmohamad
Sent: Tuesday, January 11, 2011 9:56 AM
To: vizing@paco.net
Subject: follow-up on . . . Total Chromatic Number Conjecture

Dear Professor Vizing:

Thank you very much indeed for your e-mail sent on Dec. 31st, 2010 concerning Total Chromatic Number Conjecture, TCNC. As expected, you are an excellent mathematician as well as a man of dignity and honesty. Considering the situation, neither you nor Professor Behzad who was retired officially about three decades ago need any more scientific credits. I think, with your guidelines and support, I will be able to close the TCNC file, once for all, in an honorable manner. In order to do so, I would like to ask two simple questions.

Based on what I know, I am sure you will admit that Behzad was the first who declared about total coloring. Wouldn't you?

Let us ignore the fact that Behzad has introduced or studied total graphs, total groups, total Ramsey numbers, total domination, etc. Don't you agree with me, Behzad who declared TCNC three years before you "published about it first "deserves to be considered the sole author of TCNC?

Anyhow, authors and researchers definitely should be advised to stop referring to the wrong referrals. For certain all these considerations stem from the mishandling mentioned in my e-mail dated 19 Dec. 2010 and the fact that mathematics literature should be flawless. Your comments and guidelines will be highly appreciated. With all the best wishes and regards,

Hossein Shahmohamad

# *Concluding Remarks*

I greatly admire Professor Vizing's accomplishments and deeply believe that his slip in this dispute is mere negligence.

In addition to the points mentioned in the abstract I emphasize that in spite of the negligence and the mishandling specified in this report, " Total Insanity around Total Chromatic Number Conjecture " still is an inappropriate and totally unacceptable title for Section 16.2 of The Mathematical Coloring Book [4].

In order to avoid the peril in scientific writings, young researchers and authors should be notified to read articles and essays on the serious matters such as: Plagiarism, Fabrication and Falsification.

It is our professional and academic duties to keep mathematics literature flawless. So, those who have propagated the wrong referrals [5],[6],[7], or [8] to give credit for the authorship of the TCNC to Professor V.G. Vizing are advised to compensate.

Bibliography

1: Behzad M., Graphs and their Chromatic Numbers, Doctoral Thesis, Michigan State University (1965).
2: Behzad M. and Chartrand G., Introduction to total graphs, P. Rosentiehl (ed.) Theory of Graphs, International Symposium, Rome (July 1966). Gordon and Breach, New York, 31-33 (1967).
3: G. Gutin and B. Toft, Interview with Vadim G. Vizing, European Mathematical Society (December 2000).
4: Soifer A., The Mathematical Coloring Book, Springer (2009).
5: Vizing V.G., On an estimate of the chromatic class of a p-graph, (Russian) Diskret. Analiz 3, 25-30 (1964).
6: Vizing V.G., The chromatic class of multigraph, (Russian) Kibernetika 3, 29-39 (1965).
7: Vizing V.G., Some Unsolved Problems in Graph Theory, (Russian) English Translation in Russian Math. Survey 23, 117-134 (1968).
8: Zykov A.A., Problem 12 (by V.G. Vizing), beitrage zur Graphebthorie,Vorgetragen auf dem Internationalen Kolloquium in Manebach (DDR) Vom 9-12, Mai 1967, P. 228, B.G. Teubner, (1968).

9: Zykov A.A., On Some New Results of Soviet Mathematicians, P. Rosentiehl (ed.) Theory of Graphs. International Symposium, Rome (July 1966), Gordon and Breach, New York, 415-416 (1967).